\documentclass[review]{elsarticle}
\makeatletter
\def\ps@pprintTitle{%
 \let\@oddhead\@empty
 \let\@evenhead\@empty
 \def\@oddfoot{\centerline{\thepage}}%
 \let\@evenfoot\@oddfoot}
\makeatother

\usepackage{graphicx}
\usepackage{amsfonts, amsmath, amssymb}
\usepackage{booktabs,siunitx}
\usepackage[svgnames,table]{xcolor}
\usepackage{bm}
\usepackage{cancel}
\usepackage{color}
\usepackage{enumerate}
\usepackage{float}
\usepackage{hyperref}
\usepackage[english]{babel}
\usepackage[margin=1in]{geometry}
\usepackage{mathtools}
\usepackage{setspace}
\usepackage{amsthm}
\usepackage{subfigure}
\usepackage[tableposition=above]{caption}
\usepackage{multirow}
\usepackage{lineno}
\usepackage{array,multirow}
\newcolumntype{K}[1]{>{\centering\arraybackslash}p{#1}}
\newcommand{\RN}[1]{%
	\textup{\uppercase\expandafter{\romannumeral#1}}%
}

\newcommand \D [2]{\frac{\partial #1}{\partial #2}}

\renewcommand{\vec}[1]{\bm{\mathrm{#1}}}
\newcommand{\V}[1]{\bm{\mathrm{#1}}}

\def \x{\vec{x}}

\def \cA{\vec{\mathcal{A}}}
\def \cG{\vec{\mathcal{G}}}
\def \cW{\vec{\mathcal{W}}}

\def \cP{{\mathcal{P}}}

\def \cL{{\mathcal{L}}}

\def \half{\frac{1}{2}}
\def \3half{\frac{3}{2}}

\def \x{\vec{x}}

\makeatletter
\renewcommand*\env@matrix[1][\arraystretch]{%
  \edef\arraystretch{#1}%
  \hskip -\arraycolsep
  \let\@ifnextchar\new@ifnextchar
  \array{*\c@MaxMatrixCols c}}
\makeatother

\begin{document}

\let\today\relax

\begin{frontmatter}

\title{On immersed boundary kernel functions: a constrained quadratic minimization perspective}

\author[SDSU]{Amneet Pal Singh Bhalla\corref{mycorrespondingauthor}}
\ead{asbhalla@sdsu.edu}

\address[SDSU]{Department of Mechanical Engineering, San Diego State University, San Diego, CA}
\cortext[mycorrespondingauthor]{Corresponding author}

\begin{abstract}
In the immersed boundary (IB) approach to fluid-structure interaction modeling, the coupling between the fluid and  structure variables is mediated using a regularized version of Dirac delta function. In the IB literature, the regularized delta functions, also referred to IB kernel functions, are either derived analytically from a set of postulates or computed numerically using the moving least squares (MLS) approach. Whereas the analytical derivations typically assume a regular Cartesian grid, the MLS method is a meshless technique that can be used to generate kernel functions on complex domains and unstructured meshes. In this note we take a viewpoint that IB kernel generation, either analytically or via MLS, is a constrained quadratic minimization problem. The extremization of a constrained quadratic function is a broader concept than kernel generation and there are well-established numerical optimization techniques to solve this problem. For example, we show that the constrained quadratic minimization technique can be used to generate one-sided (anisotropic) IB kernels and/or to bound their values.    
\end{abstract}


\end{frontmatter}


\section{Introduction} \label{sec_intro}

The immersed boundary (IB) method was originally proposed by Peskin to model blood flow around heart valves~\cite{Peskin1977}. In the IB formulation, the fluid occupying the entire computational domain is described by Eulerian variables, whereas  
the immersed structure that is allowed to freely cut the background Eulerian grid is described by Lagrangian variables. The presence of the immersed structure leads to an additional body force term in the momentum equation. The Lagrangian structure moves with the local fluid velocity and satisfies the no-slip boundary condition on the fluid-structure interface. The Eulerian-Lagrangian coupling in the IB formulation is mediated using regularized integral kernels. Specifically, IB forces defined on the marker points are \emph{spread} to the Eulerian grid from the Lagrangian structure, and conversely, the fluid velocity is \emph{interpolated} from the Eulerian mesh to the marker points. The use of an isotropic and regularized IB kernel \emph{diffuses} the fluid-structure interface to nearby grid cells and this smearing of the interface is proportional to the width of the kernel. This also implies that the (usual) isotropic IB kernels couple the structural degrees of freedom to fluid variables on both sides of the interface.

In the literature, IB kernels have been derived either analytically or computed numerically via the moving least squares (MLS) approach. Peskin~\cite{Peskin02} used a set of postulates to derive analytical expressions for several IB kernels such as the four-, five-, and six-point function. These postulates aim to achieve several desirable features in an IB kernel, such as polynomial reproducing conditions (also known as the discrete moment conditions), grid translational invariance property, and compact kernel support. More recently, Bao et al.~\cite{Bao2016} used a weaker second moment condition to derive improved versions of the standard five-, and six-point IB kernels. A weaker second moment condition allowed the authors to remove the negative tails of the standard five- and six-point  kernels, and to enhance their gird translational invariance property. A main limitation of the analytical approach is the assumption of a regular Cartesian grid  and the tensor-product form of the kernel function in dimensions higher than one. To circumvent these limitations, Vanella and Balaras introduced the MLS approach to generate IB kernel functions~\cite{Vanella09}. Although the  authors in~\cite{Vanella09} demonstrated the IB/MLS approach on a regular Cartesian grid, the technique is easily extendable to unstructured grid-based IB methods as shown by  Krishnan et al.~\cite{Krishnan17} and Saadat et al.~\cite{Saadat18}. In this work we take a viewpoint that IB kernel generation is a constrained quadratic
minimization problem and both analytical and MLS formulations can be described using the language of  optimization. 

In the prior IB/MLS works~\cite{Vanella09,Li15,Tullio16,Krishnan17,Le2017,Saadat18,Haji19}, a standard formulation of the moving least squares method was used to transform a cubic spline weight function into a regularized IB kernel function that satisfies zeroth- and first-order discrete moment conditions. The standard MLS formulation minimizes the weighted $L^2$-norm of the difference between the interpolant and the given data values. A lesser known formulation of the moving least squares problem (in the IB literature) is the Backus-Gilbert theory, which was developed in the context of geophysics~\cite{Backus1968resolving}. In the  Backus-Gilbert formulation, polynomial reproducing conditions are used in the MLS problem definition. However, both formulations ultimately lead to the same IB kernel function as shown in this work. The connection between the two MLS  formulations was first pointed out by Bos and \v{S}alkauskas in~\cite{Bos1989moving} and re-iterated here in the context of constrained quadratic minimization. We also show that the constrained quadratic minimization approach can be used to generate discrete moment conditions satisfying anisotropic IB kernels that effectively couple structural degrees of freedom to fluid variables on only one side of the fluid-structure interface. This can help to reduce spurious feedback forcing and ``leakage" that are typically observed in the diffuse-interface IB models. Furthermore, the constrained quadratic minimization method can be used to bound the weights of the IB kernels.   

In the following sections, we first review the two formulations of the moving least squares problem and describe them using the language of optimization. Next, we review the set of postulates used by Peskin to derive the four-point IB kernel and frame it into the language of constrained quadratic minimization. Finally, we describe methods to generate one-sided and/or bounded IB kernels and demonstrate this process using a two-dimensional example.
  

\section{Moving least squares method} \label{sec_mls_theory}

\subsection{The standard formulation} \label{sec_stdMLS}

The standard formulation of the moving least squares method considers the following \emph{approximation} problem. Given a (column) vector of data values $\V{g} = (g(\x_1), \ldots, g(\x_N))^\intercal$, defined on $N$ distinct data sites $\bm{\mathcal{X}} = \{\x_1,\ldots,\x_N\}$ by a (smooth) function $g$, find the quasi-interpolant $\cP g$ to $g$
\begin{equation}
\label{eqn_stdmls}
\cP g(\x) = \sum_{j= 1}^m c_j(\x) p_j(\x)  = \V{p}^\intercal (\x)\V{c}(\x),      
\end{equation}
at an \emph{evaluation} point $\x$ using the polynomial approximation space  $\mathbb{P} = \{p_1,\ldots, p_m\}$. Typically, $ m \ll N$ and the coefficient vector $\V{c}(\x)$ in Eq.~\eqref{eqn_stdmls} is found by minimizing the weighted $L^2$-norm of the error function $J(\V{c})$
\begin{equation}
 \label{eqn_stderror}
J(\V{c}(\x)) = \sum_{i= 1}^N  W(\x_i,\x) [\V{p}^\intercal (\x_i)\V{c}(\x) - g(\x_i)]^2,      
\end{equation}   
with respect to $\V{c}$. In the equation above, $W(\x) = W(\x_i,\x)$ is a given weight function which decreases in magnitude away from the evaluation point and is taken to be non-negative. The error  function $J$ is minimized using the \emph{essential} condition of extremum, $\D{J}{\V{c}} = \V{0}$, which leads to the normal equation of the form
\begin{equation}
\label{eqn_stdmls_normal}
 \sum _{j = 1}^m c_j(\x) \left< p_j,p_k \right>_{W(\x)} = \left<g, p_k\right>_{W(\x)}, \quad k = 1,\ldots,m.
\end{equation}
Here, $\left< \cdot, \cdot \right>_{W(\x)}$ denotes the weighted inner product of two quantities. If $\cA \in \mathbb{R}^{m \times n}$ denotes the polynomial matrix with entries $\cA_{ij} = p_i(\x_j), i = 1, \ldots,m, j = 1, \ldots, N$, and $\cW(\x) =  \text{diag}\left( \V W \right)$ is the diagonal matrix containing weights  
(with $\V{W} = \left[W(\x_1,\x), \dots, W(\x_N,\x) \right]$ as the main diagonal of $\cW$), then in matrix notation Eq.~\eqref{eqn_stdmls_normal} reads as
\begin{equation}
\label{eqn_stdmls_gram}
\cG(\x) \V{c}(\x) =  \V{g}_p(\x).
\end{equation}
Here, $\cG$ is the positive-definite Gram matrix $\cG(\x) = \cA \cW \cA^\intercal \in \mathbb{R}^{m \times m}$ and $\V{g}_p$ is the right-hand side vector $\V{g}_p(\x) = \cA \cW \V{g} \in \mathbb{R}^{m}$. Note that the coefficient vector $\V{c}(\x) = \cG^{-1}\V{g}_p$, Gram matrix $\cG(\x)$, and the right hand side vector  $\V{g}_p(\x)$, all depend on the evaluation point $\x$. Consequently, the normal Eq.~\eqref{eqn_stdmls_gram} is solved anew each time the evaluation point is changed or ``moved."  

In the context of the immersed boundary method, the standard MLS formulation is used to interpolate velocity from Eulerian grid nodes to a Lagrangian marker point and to spread IB force from a Lagrangian marker point to the surrounding Eulerian nodes~\cite{Vanella09}. This is achieved by considering a particular Lagrangian marker point as the current evaluation point and the surrounding Eulerian grid nodes (to the Lagrangian point) as the distinct data sites that provide velocity or force data. The IB/MLS kernel function  $\V{\Psi}(\x)  = \left(\psi_1(\x), \ldots,\psi_N(\x)\right)^\intercal$, also known as the \emph{generating function}, is obtained by combining Eqs.~\eqref{eqn_stdmls} and~\eqref{eqn_stdmls_normal}
\begin{align}
\cP g(\x) = \sum_{i= 1}^N \psi_i(\x) g(\x_i) = \V{g}^\intercal \V{\Psi} (\x)  &= \V{p}^\intercal (\x)\V{c}(\x)  = \V{g}^\intercal \cW(\x) \cA^\intercal(\x) \cG^{-1}(\x) \V{p}(\x).
\end{align}
Therefore, from the equation above we get 
\begin{equation}
\label{eqn_stdmls_psi}
\V{\Psi} (\x) = \cW (\x) \cA^\intercal(\x) \cG^{-1}(\x) \V{p}(\x). 
\end{equation}    
Comparing the two forms of the quasi-interpolant $\cP g$ in Eqs.~\eqref{eqn_stdmls} and~\eqref{eqn_stdmls_psi}, it is clear that the latter form employs a kernel function $\V{\Psi} = \V{\Psi}(\x)$ that depends only on the evaluation point $\x$. In contrast, the former form uses a coefficient vector $\V{c} = \cG^{-1}\V{g}_p$, which depends on both the evaluation point $\x$ and the data values $\V{g}$, i.e., $\V{c} = \V{c}(\x;\V{g})$. For IB methods, it is therefore more convenient to work with $\V{\Psi}$ directly, as it allows the same kernel function to be used in both velocity interpolation and force spreading operations. It can be shown that using the same IB kernel function in the two operators conserve energy during Lagrangian-Eulerian coupling~\cite{Peskin02}.

Although, it is not explicitly enforced in the error minimization of $J$, the standard MLS method will reproduce the  polynomial basis functions $p_1,\ldots,p_m$, i.e., 
\begin{equation}
\label{eqn_stdmls_pr}
\cP p_j(\x)  = \V{p}^\intercal (\x)\V{c}(\x; \V{p}_j) =  \V{p}_j^\intercal(\x_i) \V{\Psi} (\x)= p_j(\x), \quad \quad \forall p_j \in \mathbb{P}.      
\end{equation}
We shall prove this property of standard MLS formulation in Sec.~\ref{sec_BGMLS}.

Next, we show that the standard MLS formulation can be viewed as a constrained quadratic minimization problem of the form
\begin{equation}
\label{eqn_probA}
\text{Problem A} = \begin{cases}
  \min\limits_{\V{c}(\x) \in  \mathbb{R}^{m} }  & J  \; = \; \half\, \V{c}^\intercal(\x) \cG(\x) \V{c}(\x) \\
  \text{subject to:} & \cG(\x) \V{c}(\x) =  \V{g}_p(\x).
  \end{cases}
\end{equation} 
Note that Problem A is a minimization problem because $\cG$  is a (symmetric) positive-definite matrix. The above constrained optimization problem can be solved by introducing the Lagrangian function $\cL$, with the help of a Lagrange multiplier $\V{\lambda}(\x) \in \mathbb{R}^m$ 
\begin{equation}
\label{eqn_LA}
\cL = \half \, \V{c}^\intercal(\x) \cG(\x) \V{c}(\x) - \V{\lambda}^\intercal(\x)[\cG(\x) \V{c}(\x) - \V{g}_p(\x)].
\end{equation} 
Using the \emph{essential} conditions of extrema, $\D{\cL}{\V{c}} = \V{0}$ and $\D{\cL}{\V{\lambda}} = \V{0}$, we get the following equalities
\begin{align}
\label{eqn_dLA}
\cG(\x) \V{c}(\x) &= \cG(\x) \V{\lambda}(\x),   \\
\cG(\x) \V{c}(\x) &=  \V{g}_p(\x),
\end{align} 
respectively. Solving the above two equations, we obtain  $\V{c}(\x) = \V{\lambda} (\x) = \cG^{-1}(\x)\V{g}_p(\x)$. Hence, the solutions to the standard MLS formulation and Problem A are the same.

\subsection{The Backus-Gilbert formulation} \label{sec_BGMLS}

The Backus-Gilbert  formulation of the moving least squares method seeks the quasi-interpolant of the type
\begin{equation}
\label{eqn_bgmls}
\cP g(\x) = \sum_{i= 1}^N g(\x_i) \psi_i(\x),      
\end{equation}
that satisfies the polynomial reproducing conditions 
\begin{equation}
\sum_{i= 1}^N p(\x_i) \psi_i(\x) = p(\x), \quad \quad \forall p \in \mathbb{P}.  \label{eqn_poly_reproduce}
\end{equation}
The polynomial reproduction constraints correspond to the \emph{discrete moment conditions} for the function $\psi_i(\x)$.  Eq.~\eqref{eqn_poly_reproduce} in matrix form reads as $\cA \V{\Psi}(\x) =  \V{p}(\x)$, in which $\cA$ is the same  polynomial matrix as defined in the previous section~\ref{sec_stdMLS}. The  desired generating function $\V{\Psi}$ in the Backus-Gilbert  formulation is obtained by solving the following constrained quadratic minimization problem: 
\begin{equation}
\label{eqn_probB}
\text{Problem B} = \begin{cases}
  \min\limits_{\V{\Psi}(\x) \in  \mathbb{R}^{N} }  & J  \; = \; \half\, \V{\Psi}^\intercal(\x) \cW^{-1}(\x) \V{\Psi}(\x) \\
  \text{subject to:} & \cA \V{\Psi}(\x) =  \V{p}(\x).
  \end{cases}
\end{equation} 

The solution to Problem B is obtained by minimizing the Lagrangian function $\cL$, given by
\begin{equation}
\label{eqn_LB}
\cL = \half\, \V{\Psi}^\intercal(\x) \cW^{-1}(\x) \V{\Psi}(\x) - \V{\lambda}^\intercal(\x)[\cA \V{\Psi}(\x) -  \V{p}(\x)], 
\end{equation} 
in which $\V{\lambda}(\x) \in \mathbb{R}^m$ is the Lagrange multiplier. The Lagrangian $\cL$  is minimized by using the essential conditions of extrema, $\D{\cL}{\V{\Psi}} = \V{0}$ and $\D{\cL}{\V{\lambda}} = \V{0}$, which yields the following equalities
\begin{align}
\label{eqn_dLB_1}
\cW^{-1}  \V{\Psi} - \cA^\intercal \V{\lambda} & = \V{0}, \\
\label{eqn_dLB_2}
\cA \V{\Psi} - \V{p} &= \V{0}.
\end{align} 
respectively. Solving the above two equations, we obtain 
\begin{align}
\label{eqn_BG_lambda}
\V{\lambda} & =   \left( \cA \cW \cA^\intercal \right)^{-1} \V{p} = \cG^{-1} \V{p}, \\
\label{eqn_BG_psi}
 \V{\Psi}  &= \cW \cA^\intercal \V{\lambda} = \cW \cA^\intercal \left( \cA \cW \cA^\intercal \right)^{-1} \V{p} = \cW \cA^\intercal \cG^{-1} \V{p}.
\end{align} 

Comparing Eqs.~\eqref{eqn_BG_psi} and~\eqref{eqn_stdmls_psi}, it is clear that both formulations of the moving least squares method yield the same quasi-interpolant $\cP g$ to function $g$. Since the polynomial reproduction constraints are directly included in the quadratic function minimization used to solve Problem B, the equivalence of the two formulations implies that the standard moving least squares method also reproduces the polynomial basis functions. 

In the IB literature, constant and linear polynomials are typically employed in the moving least squares formulation. Defining the polynomials relative to the evaluation point $\x = (\bar{x},\bar{y}, \bar{z})$, the set of basis functions $\mathbb{P}$ for IB problems consists of 
\begin{equation}
\label{eqn_polybasis}
\mathbb{P} = \{1, x - \bar{x}, y- \bar{y}, z - \bar{z} \}.
\end{equation} 

Note that if the weight function $W(\x_i,\x)$ already satisfies the discrete moment conditions on the support of the IB/MLS kernel, then the constrained quadratic minimization would result in a trivial solution $\V{\Psi} = \cW$. Consequently, the product $\cA^\intercal \V{\lambda}$ in Eq.~\eqref{eqn_BG_psi} produces a vector whose entries are all one.

\section{Peskin's IB kernel} \label{sec_ib4}

If the IB simulations are performed on a regular Cartesian grid (which is often the case), then Peskin's delta functions can be used instead of the IB/MLS kernels described earlier.  Peskin's delta functions satisfy the polynomial reproducing conditions, can be expressed compactly in a piecewise analytical manner,  and furthermore, they do not explicitly depend upon a particular evaluation/marker point. The latter property is due to the fact that Peskin's delta functions are constructed in terms of relative distance between the data site and the evaluation point, i.e,  $\psi = \psi(|\x_i - \x|)$ and $\psi \ne \psi(\x_i , \x)$. 

Taking the example of widely popular four-point delta function $\psi_{4}$ from the IB literature, we list the four postulates used by Peskin for its construction below:
\begin{align*}
(i) \text{  Zeroth moment:}   &\quad \quad\sum_{i= 1}^N \psi_4(|\x_i - \x|)  = 1,  \\
(ii) \text{  Even-odd:}   &\quad \quad \sum_{i \; \text{even}} \psi_4(|\x_i - \x|)  = \sum_{i \; \text{odd}} \psi_4(|\x_i - \x|) = \half,  \\
(iii) \text{  First moment:}   &\quad \quad \sum_{i= 1}^N (\x_i - \x) \psi_4(|\x_i - \x|)  = 0, \\
(iv) \text{  Sum-of-squares:}   &\quad \quad \sum_{i= 1}^N  \psi_{4}^2(|\x_i - \x|)  = \left(\frac{3}{8}\right)^d. \quad \quad (d \text{ denotes the spatial dimension})
\end{align*} 
Note that the even-odd condition in the second postulate implies the zeroth moment condition (or the constant  reproducing condition), and therefore, the first equation is not explicitly used in deriving $\psi_4$. The first moment condition in the third postulate is the linear polynomial reproducing condition, whereas the sum-of-squares condition in the last postulate is a weaker condition on grid translational invariance.  

Similar to the IB/MLS kernel formulation, Peskin's four-point delta function can be formulated as a constrained quadratic minimization problem. The minimization problem reads as
 \begin{equation}
\label{eqn_probC}
\text{Problem C} = \begin{cases}
  \min\limits_{\V{\Psi}(\x) \in  \mathbb{R}^{N} }  & J  \; = \; \half\, \V{\Psi}^\intercal(\x) \cW^{-1}(\x) \V{\Psi}(\x) \\
  \text{subject to:} & (\x_i - \x)^\intercal \V{\Psi}(\x) =  \V{0}, \\
                             & \V{1}_{\rm even}^\intercal \V{\Psi}(\x)     = \half, \\
                             & \V{1}_{\rm odd}^\intercal \V{\Psi}(\x)     = \half, \\
                             & \V{\Psi}^\intercal(\x) \V{\Psi} = \left(\frac{3}{8}\right)^d,
  \end{cases}
\end{equation} 
in which $\V{1}_{\rm even}$ is a (column) vector containing the values one and zero at even- and odd-indexed data sites, respectively. $\V{1}_{\rm odd}$ is defined analogously. Note that Problem C has a nonlinear equality constraint, which would require a nonlinear optimizer, such as the interior point algorithm~\cite{Wachter2006} or the sequential quadratic programming technique~\cite{Nocedal2006}. 

\section{One-sided and bounded IB kernels} \label{sec_ib4}

Next, we show that the constrained quadratic minimization approach can help to generate IB kernels that effectively couple structural degrees of freedom to fluid variables on only one side of the fluid-structure interface. Such interpolation kernels have been termed ``one-sided IB kernels" by Bale et al.~\cite{Bale2021}, who have recently proposed them. The one-sided IB kernels help to reduce spurious feedback forcing and internal flows that are typically observed in IB models that use isotropic kernel functions which couple the structure to fluid degrees of freedom on both sides of the interface~\cite{Bale2021}. However, such kernels have large positive and negative weights, which can lead to numerical instability over time. To mitigate the numerical instability the authors in~\cite{Bale2021} mollified the one-sided IB kernel weights at the expense of the first moment condition. In contrast, here we show that the numerical values of  the one-sided kernel function can be bounded, e.g., made positive or constrained within a smaller range using the constrained quadratic minimization approach, while still  satisfying the first moment condition. 

\begin{figure}
	\centering
	\subfigure[Support of an isotropic IB kernel]{
		\includegraphics[scale=0.42]{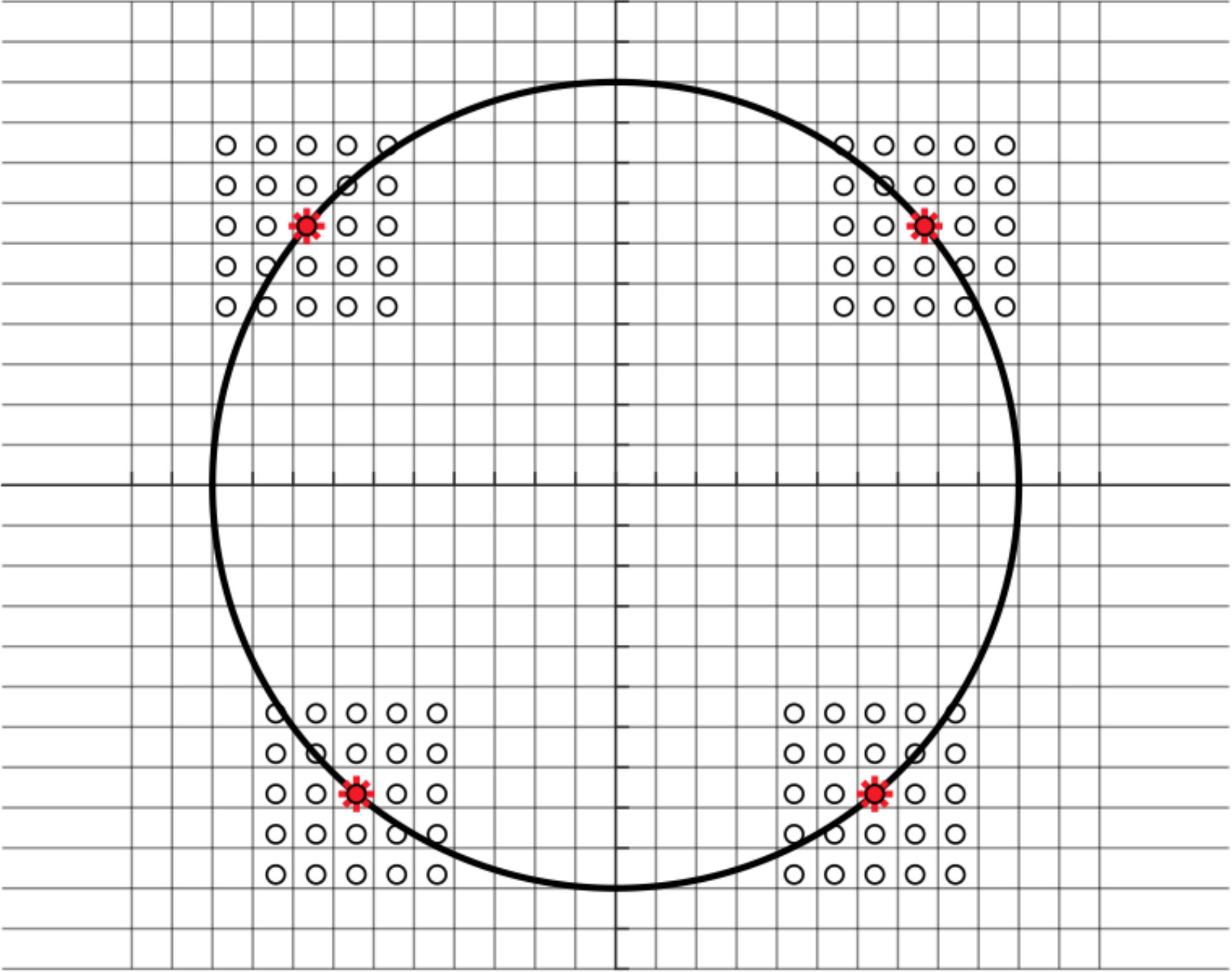} 
		\label{fig:psi_fullsupport}
	}
	\subfigure[Support of a one-sided IB kernel] {
		\includegraphics[scale=0.42]{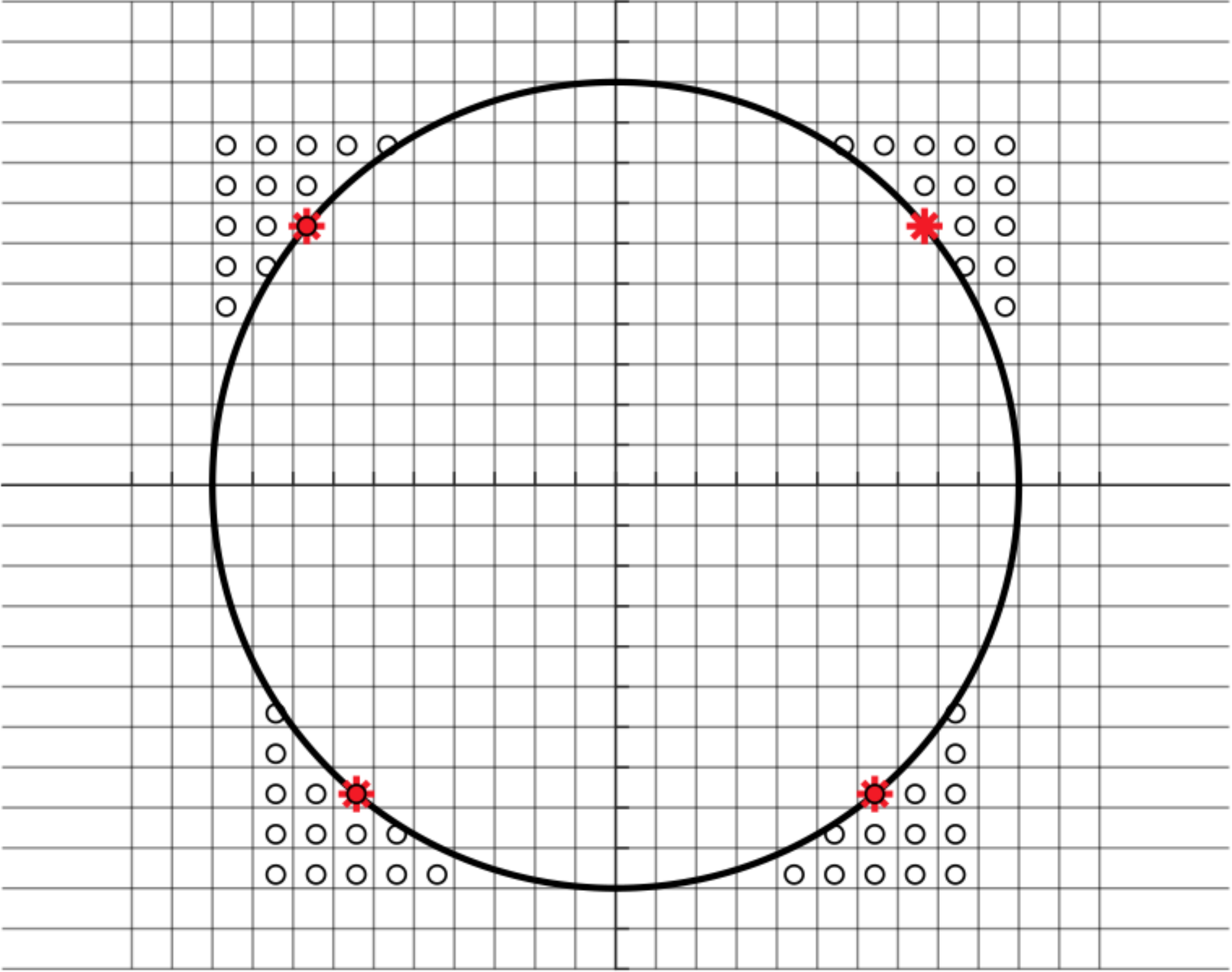} 
		\label{fig:psi_innersupport}
	}	
	\caption{Representative Lagrangian markers (${\color{red}*}$) on a circular interface embedded in a Cartesian grid.\subref{fig:psi_fullsupport} Full support of a typical isotropic IB kernel and~\subref{fig:psi_innersupport} restricted support of one-sided IB kernel. The Eulerian support of the maker points are shown by circles ($\circ$).}
	\label{fig:psi_support}
\end{figure}

We demonstrate the one-sided IB kernel construction by considering a circular interface that is embedded in a regular Cartesian grid. The interface is discretized using Lagrangian marker points, four of which are highlighted by red markers in Fig.~\ref{fig:psi_support}. Typical isotropic IB kernels consider both sides of the interface to define the support region; this is shown by black circles for the four marker points in Fig.~\ref{fig:psi_fullsupport}. For the one-sided kernel, the Eulerian support region is restricted to only one side of the interface. For example, in Fig.~\ref{fig:psi_innersupport} the kernel support is restricted to the exterior of the cylinder, which is denoted $\Omega^{+}$ region.  Analogously, the kernel support can be defined interior to the cylinder ($\Omega^{-}$ region); here we discuss the approach in the context of the former case. The one-sided, bounded, and discrete moment (polynomial reproducing) conditions satisfying IB kernels can be constructed numerically by solving a constrained quadratic minimization problem of the type
  \begin{equation}
\label{eqn_probD}
\text{Problem D} = \begin{cases}
  \min\limits_{\V{\Psi}(\x) \in  \mathbb{R}^{N} }  & J  \; = \; \half\, \V{\Psi}^\intercal(\x) \cW^{-1}_{\rm R}(\x) \V{\Psi}(\x) \\
  \text{subject to:} & \cA \V{\Psi}(\x) =  \V{p}, \\
                             & \V{\Psi}_{\rm L} \le \V{\Psi}(\x)  \le \V{\Psi}_{\rm U}  
   \end{cases}
\end{equation} 
In the Problem D above, $\cW_{\rm R}$ is the restricted version of the weight function that is defined to be
\begin{equation}
\cW_{\rm R} = \begin{cases}
W(\x_i,\x) & \text{ if } \x_i \in \Omega^{+}, \\
0               & \text{ if } \x_i \in \Omega^{-}.
\end{cases}
\end{equation} 
  
The generating function $\V{\Psi}$ in Eq.~\eqref{eqn_probD} is bounded on the lower end by $\V{\Psi}_{\rm L}$ values and on the upper end by $\V{\Psi}_{\rm U}$ values. In lack of any particular requirement, all weights are bounded alike, i.e.,  $\V{\Psi}_{\rm L} = \alpha \V{1} = \V{\alpha}$ and $\V{\Psi}_{\rm U} = \beta \V{1} = \V{\beta}$, in which $\alpha$ and $\beta$ are two scalars and $\V{1}$ is a column vector with all entries as one. Note that we have included an inequality constraint, along with the equality constraint in Problem D. Without the inequality constraint, Problem D is same as Problem A or Problem B.

Next, we solve Problem D for the circular interface example by taking  constant and linear polynomial basis functions in $\mathbb{P}$ (as in Eq.~\eqref{eqn_polybasis}). A linear field of the type
\begin{equation}
g(x,y) = 10 x + 5 y
\end{equation}
is defined on the cell centers of a uniform Cartesian mesh that discretizes the computational domain of extents $\Omega \in [-1,1]^2$. The uniform cell size is taken to be $h = 0.075$. As the circular interface partitions the domain into inner $\Omega^-$ and outer  $\Omega^+$ regions,  $\Omega = \Omega^- \cup \Omega^+$. The cylinder is centered around the point $(0,0)$ and has a radius of $R = 0.5$. The four marker points on the periphery of the cylinder make an angle of $40^\circ$, $140^\circ$ $230^\circ$, and $310^\circ$ with the $x$-axis.  A tensor-product form of the six-point spline kernel $\psi_6$ is used as a weight function in Problem D. The one-dimensional version of $\psi_6$ function reads as
\begin{equation}
\label{eqn_sixpoint}
\psi_{6}(r)  =  \begin{cases}
\frac{1}{60}( -5\kappa^5 +90\kappa^4 -630 \kappa^3 + 2130\kappa^2 - 3465\kappa+2193) & 0 \le |r| < 1, \\
\frac{1}{120}( 5\kappa^5 -120\kappa^4 +1140 \kappa^3 - 5340\kappa^2 + 12270\kappa - 10974) & 1 \le |r| < 2,  \\
\frac{1}{120}( -\kappa^5 +30\kappa^4 -360 \kappa^3 + 2160\kappa^2 - 6480\kappa+7776) & 2 \le |r| < 3,  \\
0, & 3 \le |r|.
\end{cases} 
\end{equation}
In Eq.~\eqref{eqn_sixpoint} above, $r = ( x - \bar{x})/h$ and  $ \kappa = |r| +3 $.  Note that the weight function 
$\psi_6$ satisfies constant and linear polynomial reproducing conditions when the support of the Lagrangian marker is included on both sides of the interface.  

\begin{figure}
	\centering
	\subfigure[Case 1]{
		\includegraphics[scale=0.22]{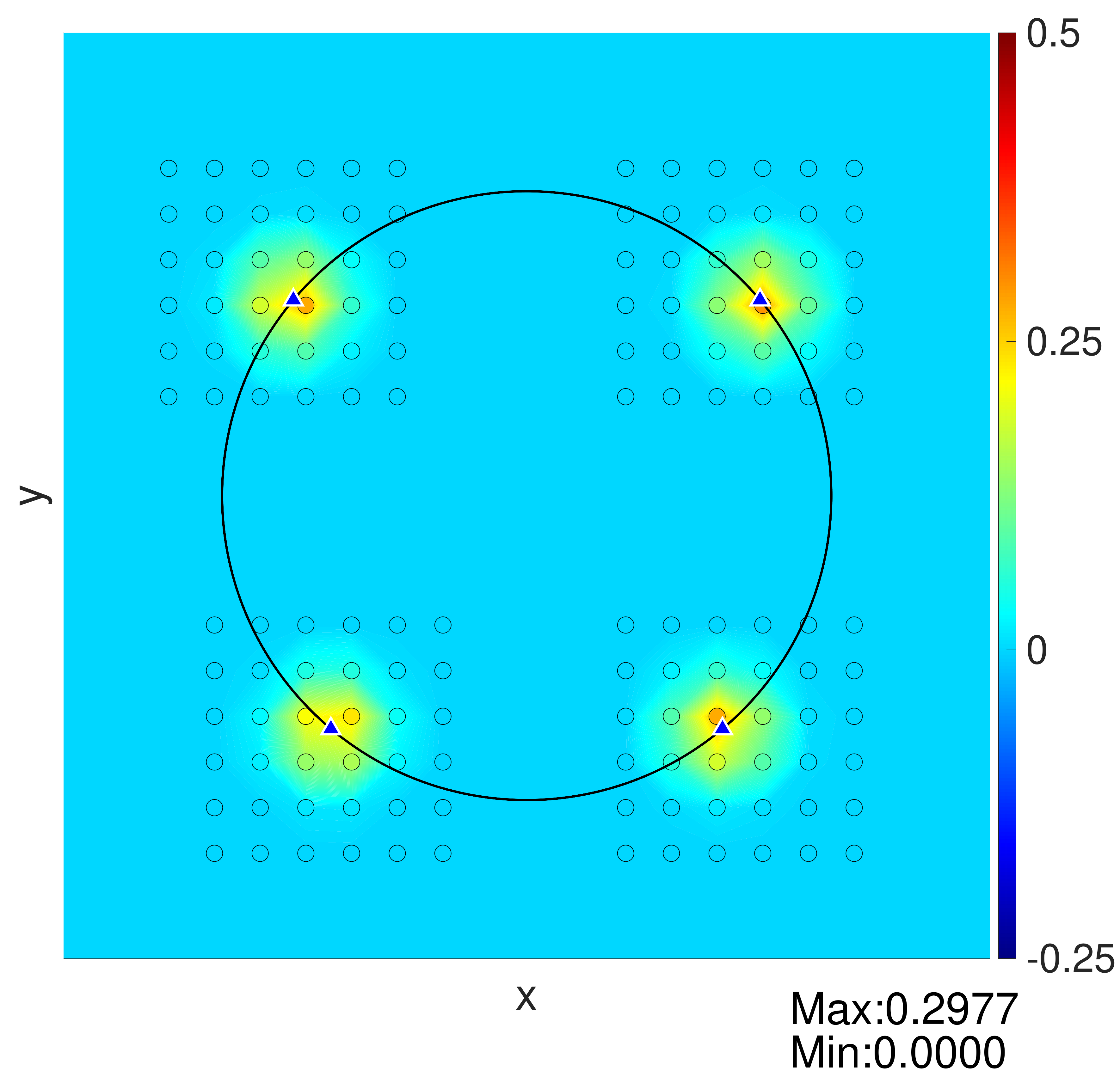} 
		\label{fig:case1}
	}
	\subfigure[Case 2] {
		\includegraphics[scale=0.22]{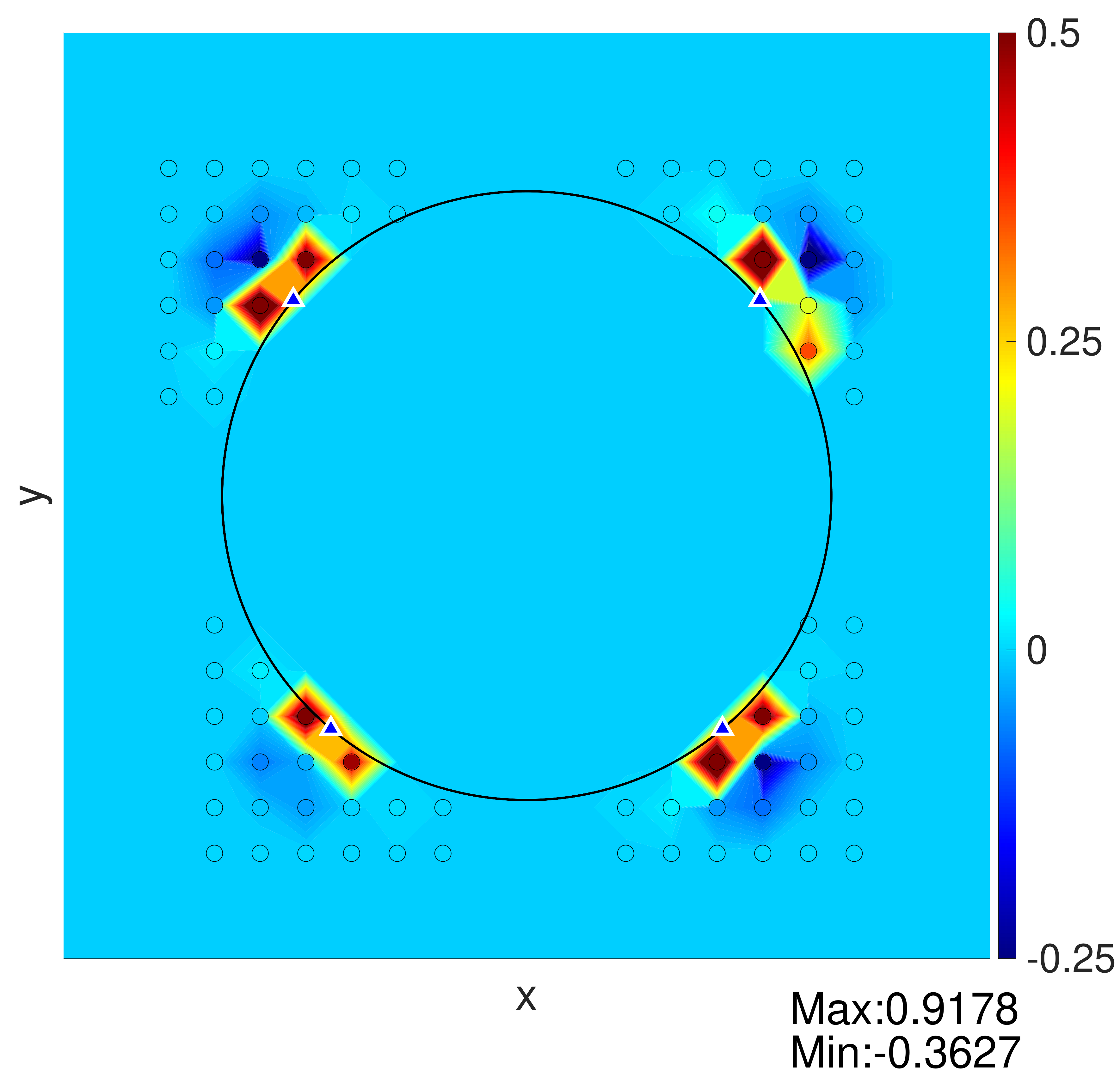} 
		\label{fig:case2}
	}	
	\subfigure[Case 3] {
		\includegraphics[scale=0.22]{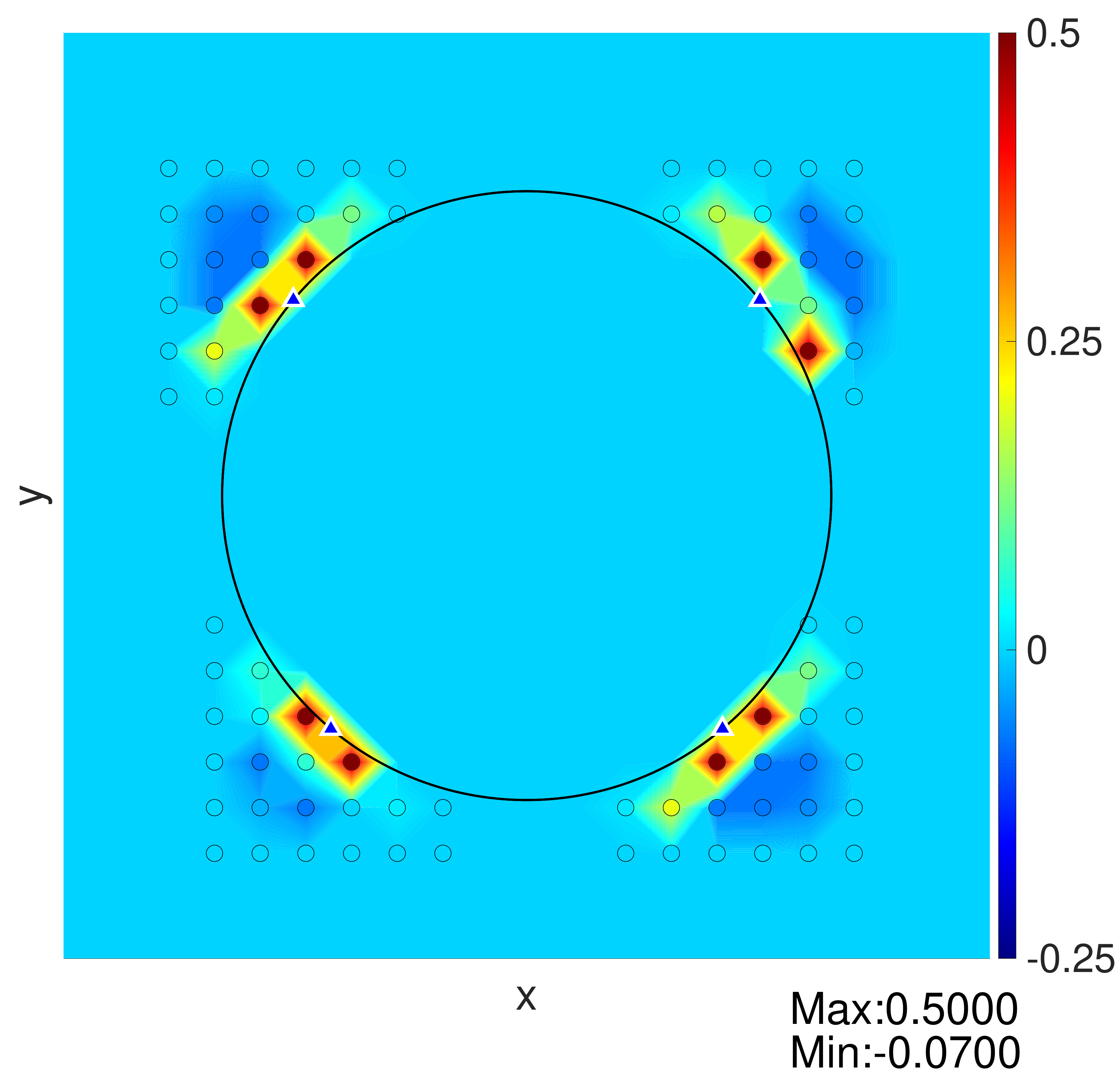} 
		\label{fig:case3}
	}
	\subfigure[Case 4] {
		\includegraphics[scale=0.22]{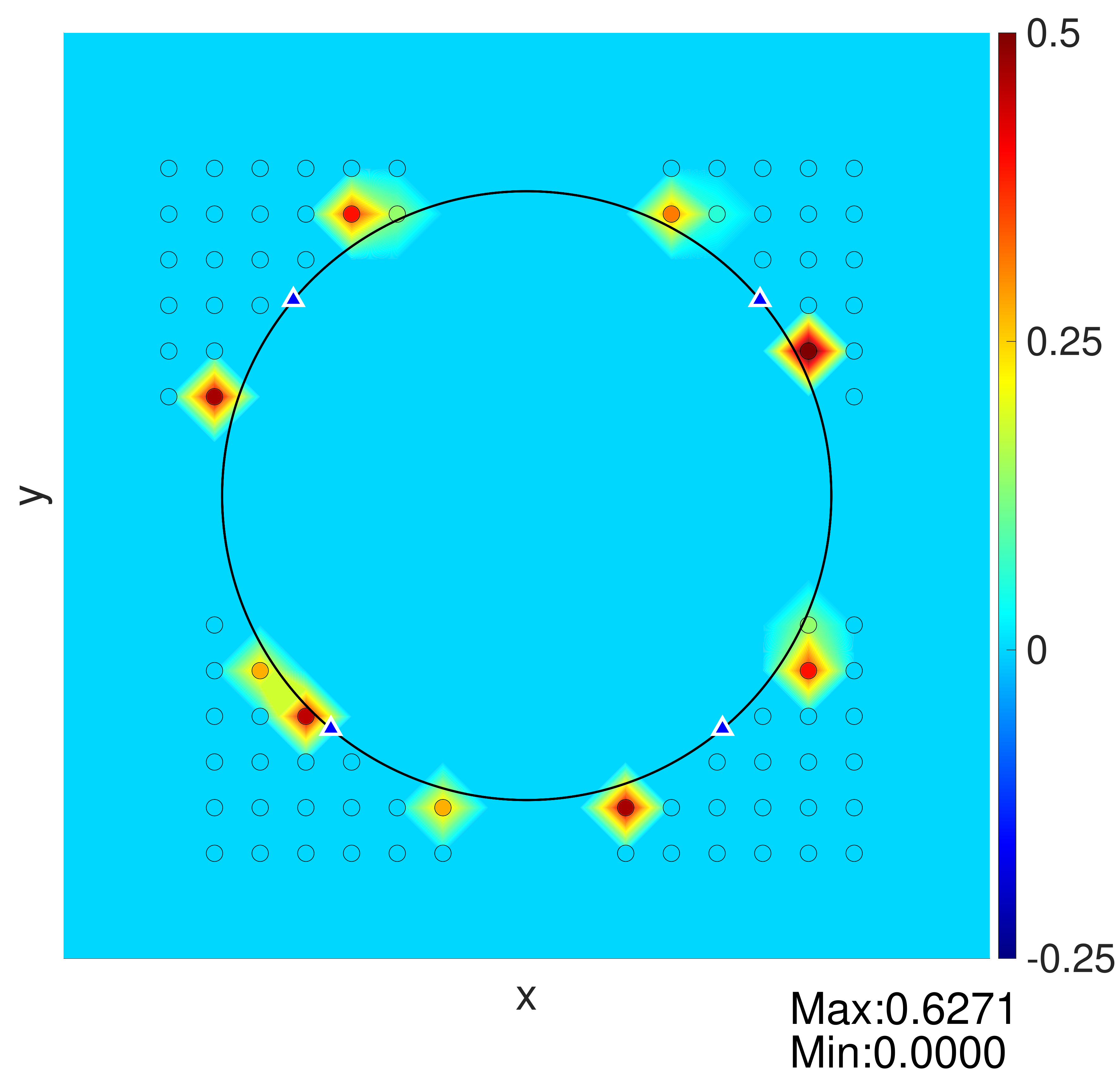} 
		\label{fig:case4}
	}	
	\caption{Color plot of $\V{\Psi}$ for four representative Lagrangian markers (\textcolor{blue}{$\blacktriangle$}) placed on a circular interface of radius $R = 0.5$. $\V{\Psi}$ is obtained using the constrained quadratic minimization formulation of Problem D. The interface is embedded in a computational domain of extents $\Omega \in [-1,1]^2$. Interpolation weights when~\subref{fig:case1} Eulerian support on both sides of the interface and only equality constraints $\cA \V{\Psi}(\x) =  \V{p}(\x) $ are used;~\subref{fig:case2} Eulerian support restricted to $\Omega^+$ region and only equality constraints $\cA \V{\Psi}(\x) =  \V{p}(\x) $ are used;~\subref{fig:case3} Eulerian support restricted to $\Omega^+$ region and both equality  $\cA \V{\Psi}(\x) =  \V{p}(\x) $ and inequality constraints $\V{-0.07}  \le \V{\Psi}(\x)  \le  \V{0.5}$ are used; and~\subref{fig:case4} Eulerian support restricted to $\Omega^+$ region and both equality  $\cA \V{\Psi}(\x) =  \V{p}(\x) $ and inequality constraints $\V{0}  \le \V{\Psi}(\x)  \le  \V{0.75}$ are used. }
	\label{fig:weights}
\end{figure}

The linear field $g$ is interpolated onto the four marker points using
\begin{itemize}

\item Case 1: Eulerian support on both sides of the interface and only equality constraints $\cA \V{\Psi}(\x) =  \V{p}(\x) $ in Problem D;
\item  Case 2: Eulerian support restricted to $\Omega^+$ region and only equality constraints $\cA \V{\Psi}(\x) =  \V{p}(\x) $ in Problem D; 
\item  Case 3: Eulerian support restricted to $\Omega^+$ region and both equality $\cA \V{\Psi}(\x) =  \V{p}(\x) $ and inequality constraints $  \V{-0.07}  \le \V{\Psi}(\x)  \le  \V{0.5}$  in Problem D; and
\item  Case 4: Eulerian support restricted to $\Omega^+$ region and both equality $\cA \V{\Psi}(\x) =  \V{p}(\x) $ and inequality constraints $  \V{0}  \le \V{\Psi}(\x)  \le  \V{0.75}$  in Problem D.
\end{itemize}
The three cases listed above are solved using the \texttt{quadprog} routine in MATLAB using the default solver options. The \texttt{quadprog} routine requires defining the cost function $J$ in terms of a matrix ($\cW^{-1}$ in this case) and lower and upper bound vectors ($\V{\Psi}_{\rm L}$ and $\V{\Psi}_{\rm R}$, respectively in this case) for the variable to be minimized ($\V{\Psi}$ in this case). 

Table~\ref{tab_error} lists the relative error in the interpolated values of $g$  and Fig.~\ref{fig:weights} shows the color plot of the interpolation weights around the four marker points. Fig.~\ref{fig:case1}, which plots the weights for Case 1 shows that Problem D generates the usual isotropic IB kernel because the support of the marker points is defined on both sides of the interface. In fact for Case 1, Problem D did not modify the underlying weight function $\V{\Psi}_6$ because the weight function already satisfies the polynomial reproducing conditions, as mentioned earlier. This can be verified by computing the $L^2$- and $L^\infty$-norms of the difference $||\V{\Psi} - \V{\Psi}_6 ||$, which come out to be $2.5829 \times 10^{-13}$ and $1.0819 \times 10^{-13}$, respectively.

Interpolation weights for Case 2 are shown in Fig.~\ref{fig:case2}. In this case we observe large positive weights near the marker points, which are followed by large negative values; $\V{\Psi}$ ranges from -0.3627 to 0.9178 in Case 2. This kind of weight distribution (having large positive and negative values) was shown to have a destabilizing effect when used within a direct-forcing IB method~\cite{Bale2021}. To make use of the one-sided IB kernel, the authors in~\cite{Bale2021} mollified the weights through a renormalization procedure. However, the mollified weights satisfy only the zeroth moment condition, i.e, they are only first-order accurate. 

Next, we show that Problem D can be used to systematically mollify the non-smooth weights obtained in Case 2 by prescribing the upper and lower bounds, without having to forfeit the first moment condition. The resulting smooth weights for Case 3 are shown in Fig.~\ref{fig:case3}, wherein we take $\alpha = -0.07$ and $\beta = 0.5$, which is a subset of the range observed for $\V{\Psi}$ in Case 2. To quantitatively access the accuracy of the IB kernels, we tabulate the interpolation error $ \mathcal{E}$ at the four marker locations in Table~\ref{tab_error}. From the table, it can be inferred that close to machine precision error norms are obtained for Case 1 and 2, which indicates that the linear field is interpolated exactly onto the Lagrangian marker points. This is not a surprise, because the linear polynomial reproducing constraint is imposed in the problem formulation itself, but it demonstrates the procedure to obtain discrete moment conditions satisfying one-sided kernels using quadratic programming, which is different compared to the moving-least squares technique of Problem A or B. For Case 3, $\mathcal{E}$ is in the range of $10^{-10} - 10^{-11}$. We remark that close to machine precision error norms can be achieved in this case by further decreasing (increasing) the lower (upper) bound of $\V{\Psi}$. However, this will not lead to any practical difference in the results of an IB simulation and loosing some digits of interpolation accuracy should be preferred over accurate but non-smooth weights, such as those shown in Fig.~\ref{fig:case2}.

Lastly,  we consider Case 4, in which our objective is to completely eliminate the negative weights of Case 2. To achieve this, we take $\alpha = 0$ and $\beta = 0.75$. Table~\ref{tab_error} shows the error norms for this case. Here, we observe that although the error norms are small (in the range of $10^{-7} - 10^{-10}$), they are not close to machine precision, indicating that the linear field is not interpolated exactly. This is because the \texttt{quadprog} routine is unable to find weights that can satisfy the equality constraints exactly, while simultaneously respecting the imposed bounds; changing the solver settings did not affect the final convergence of \texttt{quadprog} (data not shown). On the other hand, if we observe the min/max limits of the color bar in Fig.~\ref{fig:case4}, $\V{\Psi}$ for Case 4 stays within the imposed limits.  Note that although the one-sided weights of Case 4 are positive, they are not smooth and are highly localized near the interface. This type of weight distribution can also lead to numerical instability for the direct-forcing IB method~\cite{Uhlmann2005}.  Apparently, the weight-shifting technique of Bale et al. to completely eliminate the negative weights at the expense of first moment condition is more suitable for the direct-forcing IB approach than the inequality constraint employed in Case 4. Nevertheless, if we allow small magnitude of negative weights in the constrained quadratic formulation, we can achieve smoothness and as well as the desired accuracy in the one-sided IB kernels (Case 3).

\begin{table}[]
 \centering
 \caption{Relative error $ \mathcal{E} = \frac{\left| \mathcal{P}g -g \right|}{|g|}$ for the three cases of Problem D. For Case 3, $\alpha = -0.07$ and $\beta = 0.5$ and for Case 4, $\alpha = 0$ and $\beta = 0.75$}
\begin{tabular}{c c c c c }
\toprule
\multirow{2}{*}{Case}   & \multicolumn{4}{c}{Relative error $(\mathcal{E})$ at marker location}    \\
\cline{2-5}
                                    &  40$^\circ$  &  140$^\circ$  &  230$^\circ$  &  310$^\circ$      \\
\midrule
                             
    1                        & $1.6335 \times 10^{-16}$ & $1.9975 \times 10^{-16}$ & $6.9267 \times 10^{-16}$ & $1.1967 \times 10^{-15}$ \\   
\cellcolor{gray!10} 2 & \cellcolor{gray!10} $7.6776 \times 10^{-15}$ &  \cellcolor{gray!10} $2.3970 \times 10^{-14}$ &  \cellcolor{gray!10} $3.117 \times 10^{-15}$ &  \cellcolor{gray!10} $5.2997 \times 10^{-15}$    \\
  3                        & $6.1497 \times 10^{-11}$ & $1.3954 \times 10^{-10}$ & $2.2942 \times 10^{-11}$ & $4.8383 \times 10^{-10}$ \\  
 
 \cellcolor{gray!10} 4                        & \cellcolor{gray!10} $8.2015 \times 10^{-11}$ & \cellcolor{gray!10} $2.4047 \times 10^{-9}$ & \cellcolor{gray!10} $1.2423 \times 10^{-7}$ & \cellcolor{gray!10} $3.0356 \times 10^{-9}$ \\
 \bottomrule                             
 \end{tabular}
  \label{tab_error}
\end{table}


\section{Summary and Conclusions}

In this note we took a viewpoint that IB kernel generation, either analytically or via MLS, is a constrained quadratic minimization problem. The extremization of a constrained quadratic function is a broader concept than kernel generation and there are well-established numerical optimization techniques to solve this problem. We also reviewed the two formulations of the moving least squares problem and described them using the language of optimization. Further, the set of postulates used by Peskin to derive the four-point IB kernel was also described using this language. Finally, we described methods to generate one-sided and/or bounded IB kernels and demonstrated the process using a two-dimensional example. Four constrained quadratic minimization problems were described, namely, the Problem A, B, C, and D. All but Problem C involved linear constraints and can be solved using efficient optimization techniques such as the quadratic programming technique. For Problem C,  a nonlinear optimizer, such as the interior point algorithm or the sequential quadratic programming technique is required. 


\section*{Acknowledgements}
A.P.S.B~acknowledges helpful discussion with Rahul Bale. 

\appendix
\renewcommand\thesection{\Alph{section}}

\section*{Bibliography}
\begin{flushleft}
 \bibliography{References}
\end{flushleft}

\end{document}